\documentclass{amsart}
\usepackage{graphicx}
\theoremstyle{definition}

\theoremstyle{remark}

\numberwithin{equation}{section}

\begin{document}

\title{Local derivative estimates for heat equations on Riemannian manifolds}%
\author{Hong Huang}%
\address{Department of Mathematics,Beijing Normal University,Beijing 100875, P. R. China}%
\email{hhuang@bnu.edu.cn}%

\thanks{Partially supported by NSFC no.10671018.}%

\begin{abstract}
In this short note we present local derivative estimates for heat
equations on Riemannian manifolds following the line of W.-X. Shi.
As an application we generalize a second derivative estimate of R.
Hamilton for heat equations on compact manifolds to noncompact case.
\end{abstract} \maketitle

\section{Introduction}

In [S]  W.-X. Shi got  local derivative estimates for Hamilton's
Ricci flow which is very important for later developments of
geometric evolution equations, see for example Hamilton [H2] and
Perelman [P1] [P2].

In this short note we  present local derivative estimates of heat
equations on Riemannian manifolds (with static metrics) following
the line of Shi.

\hspace *{0.4cm}

{\bf Theorem 1 }  Let $\mathcal{M}^n$ be a Riemannian manifold.
Suppose $u$ is a smooth solution to the heat equation on
$B_p(2r)\times [0,T]$ for some $p \in \mathcal{M}^n$ satisfying
$|u|\leq M$. Then there exist  constants $C_k$ depending only on the
dimension and the bounds of  curvature and the covariant derivatives
( up to order $k-1$ ) of the curvature on $B_p(2r)$ such that

$|\nabla^k u|^2 \leq C_k M^2(\frac{1}{r^{2k}}+\frac{1}{t^k}+K^k),$

hold on $B_p(r)\times (0,T]$, where $K$ is the bound of curvature on
$B_p(2r)$.

\hspace *{0.4cm}

For a  similar but not  identical estimate for harmonic map heat
flow see Grayson and Hamilton [GH].

As an application we generalize a second derivative  estimate of R.
Hamilton [H1] for heat equations on compact manifolds to noncompact
case.

\hspace *{0.4cm}

{\bf Theorem 2}  Let $\mathcal{M}^n$ be a complete  Riemannian
manifold with bounded curvature and covariant derivative ( of first
order ) of Ricci curvature. Suppose  $u$ is a smooth positive
solution to the heat equation on $\mathcal{M}^n$ satisfying $u\leq
M$. Then for $0\leq t \leq 1$ we have

$t\Delta u \leq Cu[1+log(M/u)],$

where $C$ depends only on the bounds of curvature and covariant
derivative ( of first order ) of Ricci curvature.

\hspace *{0.4cm}

In the remaining two sections we prove these two theorems
respectively.

\section{Proof of  Theorem 1}

We follow the line of Shi [S] ( in particular see the exposition in
Cao-Zhu [CZ] which follows Hamilton [H3] in turn ). It is a
Bernstein-type estimate coupled with a cutoff argument.

First note that for a solution $u$ to the heat equation we have

$ \frac{\partial }{\partial t}\nabla^k u= \Delta \nabla^k u
+{\Sigma_{i=0}}^{k-1}\nabla^i Rm * \nabla^{k-i}u,$ ( in fact, when
$k=2$, we have  $ \frac{\partial }{\partial t}\nabla^2 u= \Delta
\nabla^2 u + Rm * \nabla^2 u+\nabla Rc * \nabla u )$

and

$ \frac{\partial }{\partial t}|\nabla^{k} u|^2= \Delta |\nabla^{k}
u|^2-2|\nabla^{k+1} u|^2 +{\Sigma_{i=0}}^{k-1}\nabla^i Rm *
\nabla^{k-i}u*\nabla^{k}u$

 $ \leq  \Delta |\nabla^{k} u|^2-2|\nabla^{k+1}
u|^2+C|\nabla^{k}u|^2 +C{\Sigma_{i=1}}^{k-1}|\nabla^{k-i}u|^2,$

where and below $C$ denotes various constants depending only on the
dimension and the bounds of  curvature and the covariant derivatives
of the curvature on $B_p(2r)$, and where * denotes some linear
tensor contraction, possible including constants.

Now we prove Theorem 1 by induction. Note that without loss of
generality we may assume $r\leq 1/\sqrt{K}$. We first consider the
case $k=1$.

Let $S_1(x,t)=(7M^2+u^2)|\nabla u|^2.$ Then  by Cauchy inequality we
have

$ (\frac{\partial }{\partial t}-\Delta)S_1 \leq  -|\nabla
u|^4+16M^2K|\nabla u|^2$

$\leq  -\frac{1}{2}|\nabla u|^4+ 128M^4K^2$

$ \leq -\frac{1}{64M^4}{S_1}^2+ 128M^4K^2$.

Let $F_1=\frac{S_1}{128M^4} .$  Then

$\frac{\partial F_1}{\partial t} \leq \Delta F_1 -{F_1}^2+ K^2.$

Fix a point $q \in B_p(r)$. As in [CZ] we choose a cutoff function
$\varphi$ with  support in the ball $B_q(r)$ such that $\varphi
(q)=r, 0\leq \varphi \leq Ar$ and $|\nabla \varphi|\leq A, |\nabla^2
\varphi| \leq \frac{A}{r}$, $A$ depending only on the dimension of
$M$. ( For our purpose we may pretend that $\varphi$ be smooth
everywhere by Calabi's trick.)

Let $H_1=\frac{(12+4\sqrt{n})A^2}{\varphi^2}+\frac{1}{t}+K.$  Then

$\frac{\partial H_1}{\partial t} > \Delta H_1 -{H_1}^2+ K^2.$

By maximum principle $F_1\leq H_1$ which implies $|\nabla u|^2 \leq
C_1M^2(\frac{1}{r^2}+\frac{1}{t}+K).$

Now suppose we have the bounds

$|\nabla^i u|^2 \leq C_iM^2(\frac{1}{r^{2i}}+\frac{1}{t^i}+K^i)$
   for $i\leq k$.

Then  $\frac{\partial }{\partial t}|\nabla^{k} u|^2 \leq  \Delta
|\nabla^{k} u|^2-2|\nabla^{k+1}
u|^2+CM^2(\frac{1}{r^{2k}}+\frac{1}{t^k}+K^k),$

$\frac{\partial }{\partial t}|\nabla^{k+1} u|^2 \leq  \Delta
|\nabla^{k+1} u|^2-2|\nabla^{k+2} u|^2+C|\nabla^{k+1}
u|^2+CM^2(\frac{1}{r^{2k}}+\frac{1}{t^k}+K^k).$

Let $S_k(x,t)=[B_kM^2(\frac{1}{r^{2k}}+\frac{1}{t^k}+K^k)+|\nabla^k
u|^2]\cdot|\nabla^{k+1} u|^2$. Choosing $B_k$ large enough and using
Cauchy inequality, we have

$ (\frac{\partial}{\partial t}-\Delta)S_k \leq
[-kB_kM^2t^{-k-1}-2|\nabla^{k+1}
u|^2+CM^2(\frac{1}{r^{2k}}+\frac{1}{t^k}+K^k)]\cdot |\nabla^{k+1}
u|^2+8|\nabla^k u|\cdot|\nabla^{k+1} u|^2\cdot|\nabla^{k+2}
u|+[B_kM^2(\frac{1}{r^{2k}}+\frac{1}{t^k}+K^k)+|\nabla^k
u|^2]\cdot[-2|\nabla^{k+2} u|^2+C|\nabla^{k+1}
u|^2+CM^2(\frac{1}{r^{2k}}+\frac{1}{t^k}+K^k)]$

$\leq  -|\nabla^{k+1} u|^4
+C{B_k}^2M^4(\frac{1}{r^{4k}}+\frac{1}{t^{2k}}+K^{2k})$

$\leq
-\frac{{S_k}^2}{(B_k+C_k)^2M^4(\frac{1}{r^{2k}}+\frac{1}{t^k}+K^k)^2}
+C{B_k}^2M^4(\frac{1}{r^{4k}}+\frac{1}{t^{2k}}+K^{2k})$

Let $v=\frac{1}{r^2}+\frac{1}{t}+K$ and set $F_k=bS_k/v^k .$ Then

$\frac{\partial F_k}{\partial t} \leq \Delta F_k
-\frac{{F_k}^2}{b((B_k+C_k)^2M^4v^k}+ bC{B_k}^2M^4v^k+kF_kv$

$\leq \Delta F_k -\frac{{F_k}^2}{2b((B_k+C_k)^2M^4v^k}+
\frac{1}{2}b(2C+k^2){(B_k+C_k)}^2M^4v^{k+2}.$

By choosing $b\leq 2/((2C+k^2){(B_k+C_k)}^2M^4)$, we get

$\frac{\partial F_k}{\partial t} \leq \Delta F_k
-\frac{1}{v^k}{F_k}^2+v^{k+2}.$

As in [CZ] we introduce

$H_k=5(k+1)(2(k+1)+1+\sqrt{n})A^2\varphi
^{-2(k+1)}+Lt^{-(k+1)}+K^{k+1},$

where $L\geq k+2$. Then we easily check

$\frac{\partial H_k}{\partial t} > \Delta H_k
-\frac{1}{v^k}{F_k}^2+v^{k+2}.$

By maximum principle we have

$F_k \leq H_k,$

from which one immediately get the desired estimate

$|\nabla^{k+1} u|^2 \leq C_{k+1}
M^2(\frac{1}{r^{2(k+1)}}+\frac{1}{t^{k+1}}+K^{k+1}).$

\hspace *{0.4cm}

{\bf Corollary}   Let $\mathcal{M}^n$ be a complete  Riemannian
manifold with bounded curvature and covariant derivative ( of first
order ) of Ricci curvature. Suppose  $u$ is a smooth  solution to
the heat equation on $\mathcal{M}^n$ for $t \in [0,T]$ satisfying
$u\leq M$. Then for $t \in (0,T]$ we have

$|\nabla^2 u|^2 \leq C_2 M^2(\frac{1}{t^2}+K^2),$

where $C_2$ depends only on the bounds of curvature and covariant
derivative ( of first order ) of Ricci curvature.

\section{Proof of  Theorem 2}

 The idea is to use Ni-Tam's generalized maximum principle of noncompact
 manifold in [NT]
( which  is originally  due to Karp and Li).

Let $h=\varphi[\Delta u+\frac{|\nabla u|^2}{u}]-u[n+4log(M/u)]$,
where $\varphi=(e^{Kt}-1)/Ke^{Kt}.$ As in [H1], we have

$\frac{\partial h}{\partial t}\leq \Delta h$ whenever $h\geq
0$.

It is easy to see that $\varphi \leq t$ for $t\geq 0$. Then using
Theorem 1 ( actually the Corollary ), for any  $p \in \mathcal{M}^n$
and $r>0$ we get

${\int_0}^1\int_{B_p(r)}e^{-d^2(x,p)}(\varphi\Delta u)^2dVdt$

$\leq
{\int_0}^1\int_{B_p(r)}e^{-d^2(x,p)}t^2nC_2M^2(\frac{1}{t^2}+K^2)dVdt
$

$\leq nC_2M^2(1+K^2){\int_0}^1\int_{\mathcal{M}^n}e^{-d^2(x,p)}dVdt$

$\leq  C <\infty,$

where the constant $C$ does not depend on $r$. So we get

${\int_0}^1\int_{\mathcal{M}^n}e^{-d^2(x,p)}(\varphi\Delta u)^2dVdt
< \infty.$

Combining with Kotschwar's estimate in [K] ( or one can use the case
$k=1$ of our Theorem 1 instead )  we get that

${\int_0}^1\int_{\mathcal{M}^n}e^{-d^2(x,p)}{h_+}^2dVdt < \infty,$

where $h_+(x,t):=max\{h(x,t),0\}$.

 So by the maximum principle in
[NT] we get that $h\leq 0$, and the desired  result follows.

\hspace *{0.4cm}

{\bf Remark } A similar argument was used by Kotschwar [K] to
generalize the gradient estimate of Hamilton [H1] to noncompact
case.

\hspace *{0.4cm}

{\bf Acknowledgements} I would like to thank Prof. Hongzhu Gao for
his support.

\bibliographystyle{amsplain}

\hspace *{0.4cm}

{\bf Reference}

\bibliography{1}[CZ] H.-D. Cao, X.-P. Zhu, A complete proof of the
Poincare conjecture and geometrization conjectures-application of
the Hamilton-Perelman theory of the Ricci flow, Asian J. Math. 10
(2006), 165-492.

\bibliography{2}[GH] M. Grayson, R. S. Hamilton, The formation of
singularities in the harmonic map  heat  flow, Comm. Anal. Geom. 4
(1996),no.4,525-546.

\bibliography{3}[H1] R. S. Hamilton,A matrix Harnack estimate
for the heat equation, Comm. Anal. Geom. 1 (1994),no.1,113-126.

\bibliography{4}[H2]R. S. Hamilton, A compactness property for
solution of the Ricci flow, Amer. J. Math. 117 (1995), 545-572.

\bibliography{5}[H3]R. S. Hamilton,The formation of singularities
in the Ricci flow,  Surveys in Differential Geometry 2, 7-136,
Internaional Press, 1995.

\bibliography{6}[K] B. Kotschwar,Hamilton's gradient estimate for
the heat kernel on complete manifolds, arXiv:math.AP/0701335 (to
appear in Proc. Amer. Math.  Soc.).

\bibliography{7}[NT] L. Ni,L.-F. Tam, K$\ddot{a}$hler-Ricci flow and
the Poicare-Lelong equation,  Comm. Anal. Geom. 12
(2004),no.1-2,111-141.

\bibliography{8}[P1] G. Perelman, The entropy formala for the Ricci
flow and its geometric applications, arXiv:math.DG/0211159.

\bibliography{9}[P2] G. Perelman, Ricci flow with surgery on three
manifolds, arXiv:math.DG/0303109.

\bibliography{10}[S]W.-X. Shi, Deforming the metric on complete Riemannian manifold,
JDG 30 (1989),223-301.

\end{document}